\documentclass[11pt]{amsart}

\usepackage{amsmath, amssymb, amsthm}
\usepackage{mathtools}
\usepackage{enumitem}
\usepackage{microtype}
\usepackage{hyperref}
\usepackage[capitalize,noabbrev]{cleveref}

\usepackage[backend=biber,style=alphabetic]{biblatex}
\theoremstyle{plain} 
\newtheorem{theorem}{Theorem}

\theoremstyle{definition} 

\theoremstyle{remark} 

\newcommand{\inv}[1]{#1^{-1}}
\newcommand{\nos}{\mathrm{NormalSubgroups}(G)}
\newcommand{\aand}{\quad {\textrm{and}} \quad}
\newcommand{\lean}[1]{{\small\texttt{#1}}}
\newcommand{\leandocs}{https://adomani.github.io/FittingsTheorem/docs/FittingsTheorem/NilpotentAdd.html}

\title{Fitting's Theorem and Semirings of Normal Subgroups}
\author{Damiano Testa}
\address{Mathematics Institute, University of Warwick, Coventry CV4 7AL, UK}
\email{d.testa@warwick.ac.uk}

\begin{document}

\begin{abstract}
We define a non-unital, generally non-associative, commutative semiring structure on the collection of normal subgroups of a group~$G$.
This viewpoint allows us to recast in ring-theoretic terms Fitting's classical theorem that the join of two nilpotent normal subgroups is nilpotent.
From this perspective, the two key inputs are a binomial expansion in a non-associative setting and the fact that the commutator subgroup of two normal subgroups lies in each factor.
The development is formalized in Lean, making essential use of \lean{Mathlib} for the core definitions and results.
\end{abstract}

\maketitle

\section*{Introduction}
Fitting's theorem is a foundational result in group theory.
We provide what we believe is a new approach to its proof, by recasting it in ring-theoretic terms.
This perspective allows us to unify the notion of nilpotence in a group and in a ring: we deduce Fitting's theorem from an analogous nilpotence result for sums of nilpotent elements in a ring.
Since the ring structure in question lacks some common properties (associativity, unitality and the existence of additive inverses), the formal verification proved essential in giving us confidence in this strategy.
In order to state the theorem, we recall some basic definitions.

Let~$G$ be a group and let~$H$ and~$K$ be subgroups of~$G$.
We denote the order relation of being a subgroup by $H \le G$.

We denote by $H \sqcup K$ the subgroup generated by~$H$ and~$K$.
The notation emphasizes that $H \sqcup K$ is the supremum of~$H$ and~$K$ in the lattice of subgroups of~$G$.

The subgroup~$H$ is {\emph{normal}} in~$G$ if, for every $g \in G$, the subset
\[
  gH\inv{g} = \{g h \inv{g} \mid h \in H\}
\]
coincides with $H$.

The {\emph{subgroup commutator}} of~$H$ and~$K$, denoted $[H, K]$, is the subgroup generated by all the {\emph{commutators}} of an element of~$H$ and an element of~$K$.
This means that $[H, K]$ is generated by all the elements of the form $[h, k] = h k \inv{h} \inv{k}$, with $h \in H$ and $k \in K$.

The {\emph{lower central series}} of the group~$G$ is the non-increasing sequence of groups $G = G_{0} \ge G_{1} \ge \cdots \ge G_{n} \ge \cdots$ defined by
\[
  G_0 = G \aand G_{i+1} = [G_i, G], \quad {\textrm{for }} i \ge 0.
\]
In the inductive stage, we take the commutator with the {\emph{whole group}}:
this is different from what happens with the derived series, where we would use $[G_i, G_i]$ in the recursion.
We do not consider the derived series here.

Finally, the group~$G$ is {\emph{nilpotent}} if there is a value of~$n$ such that $G_n = \{1\}$.
The smallest such~$n$, if it exists, is the {\emph{nilpotency class}} (or simply {\emph{class}}) of~$G$.

\begin{theorem}[Fitting] \label{thm:fitting}
If~$H$ and~$K$ are nilpotent normal subgroups of a group~$G$, then $H \sqcup K$ is nilpotent.
Moreover, if~$H$ has nilpotency class~$h$ and~$K$ has class~$k$, then $H \sqcup K$ has class at most $h+k$.
\end{theorem}

This paper accompanies the formalization in Lean~\cite{lean4}, available in the GitHub repository \href{https://github.com/adomani/FittingsTheorem}{FittingsTheorem}.
We equip the normal subgroups of a group~$G$ with a non-unital, generally non-associative, commutative semiring structure, with addition given by the join and multiplication given by the subgroup commutator.
Fitting's theorem then follows from a binomial theorem in this algebraic setting and standard monotonicity properties of commutators.

\section{A Semiring of Normal Subgroups}\label{sec:semiring}
Let $\nos$ denote the set of normal subgroups of a group~$G$.
We define operations
\[
H + K := H \sqcup K \qquad \text{and} \qquad H\,K := [H,K],
\]
where $[H,K]$ denotes the (subgroup) commutator.
We set $0 := \{1\}$, the trivial subgroup.
We then verify in Lean that these satisfy the axioms of a non-unital, generally non-associative, commutative semiring.
We provide an instance of \lean{NonUnitalNonAssocCommSemiring} on $\nos$.
In particular, the distributivity of multiplication over addition is a consequence of the fact that we are restricting our attention to {\emph{normal}} subgroups:
without this assumption, distributivity need not hold.
All remaining axioms of \lean{NonUnitalNonAssocCommSemiring} hold for general subgroups.

To the best of our knowledge, this semiring structure on the lattice of normal subgroups does not appear explicitly in the classical literature.

\subsection{Powers and the Lower Central Series}
A key observation is that positive powers in $\nos$ agree with terms of the lower central series:
\[
  H^k = \mathrm{lowerCentralSeries}(H, k-1) \quad (k\ge 1).
\]
This identity allows us to turn bounds for the lower central series into semiring inequalities.

\subsection{A Binomial Expansion}
We formalize a non-associative version of the binomial identity.
For this, we introduce a type of ``\lean{factor}s'' representing products of length~$n$ built by iteratively multiplying either~$H$ or~$K$ in each slot.
Given a subset $s\subseteq \{0,\dots,n-1\}$, the product $\mathrm{factor}(H,K,n,s)$ chooses~$H$ as the $i$-th factor if $i\in s$ and~$K$ otherwise.

With this notion of \lean{factor}, the formula for the power of a binomial takes the form:
\[
  (H+K)^n = \sum_{s \subseteq \{0,\dots,n-1\}} \mathrm{factor}(H,K,n,s),
\]
that is, we simply sum all possible \lean{factor}s of length~$n$ in the two letters $H,K$.
The formalization guarantees that associativity of multiplication is not needed for this formula.

\subsection{Inequalities}
The next ingredient involves inequalities.
The binomial formula produces a family of subgroups whose join contains the $(n-1)$st term of the lower central series of $H\sqcup K$.
We identify conditions implying that each \lean{factor} is trivial.
We exploit that, if~$H$ and~$K$ are normal, then $[H,K]\le H$ and $[H,K]\le K$.
By induction on~$n$, if $|s|=a$ and $|s^c|=b$ (with $a+b=n$), then
\[
  \mathrm{factor}(H,K,n,s)\le H^a
  \aand
  \mathrm{factor}(H,K,n,s)\le K^b.
\]
Thus, if the groups $H^a$ and $K^b$ are both trivial, then every term in the expansion of $(H+K)^{a+b-1}$ is trivial.

\subsection{Fitting's Theorem}
Finally, we put together all that we proved.
Suppose that~$H$ is nilpotent of class~$h$, so $H^{h+1}=\{1\}$, and similarly $K^{k+1}=\{1\}$.
Then each summand in $(H+K)^{h+k+1}$ vanishes, so $(H+K)^{h+k+1}=0$.
Consequently, the nilpotency class of $H\sqcup K$ is at most $h+k$, proving Fitting's Theorem (\cref{thm:fitting}).

\section{Related Work}
Background on Fitting's theorem can be found in standard references and surveys.
Classical statements and proofs include Fitting's original paper~\cite{Fitting1938} and modern expositions (e.g., \cite[Corollary~1.29]{Isaacs}, \cite{Hall1959, Robinson1996}).
For commutators of subgroups and their basic properties (including symmetry and monotonicity) we rely on \lean{Mathlib}'s formalization, cf.\ \cite{mathlib2020}.

\section{Artifact and Reproducibility}
The formalization is developed in Lean~4 and uses \lean{Mathlib}.
The source is available at \url{https://github.com/adomani/FittingsTheorem}.
The Lean toolchain is pinned in the file \lean{lean-toolchain} and the \lean{Mathlib} revision in \lean{lake-manifest.json}, so that the build is reproducible; the README records the build instructions.

The formalization is contained in the single file \lean{NilpotentAdd.lean}.
The statements of this paper correspond to the following declarations, each linked to its entry in the \href{\leandocs}{documentation}.
\begin{itemize}[leftmargin=2em,before=\raggedright]
  \item The semiring $\nos$ of \cref{sec:semiring}, with its instance of \lean{NonUnitalNonAssocCommSemiring}:
    \href{\leandocs\#Subgroup.NormalSubgroups}{\lean{Subgroup.NormalSubgroups}}
  \item The identification of positive powers with the lower central series:
    \href{\leandocs\#Subgroup.NormalSubgroups.pow_eq_lowerCentralSeries_pred}{\lean{Subgroup.NormalSubgroups.pow\_eq\_lowerCentralSeries\_pred}}
  \item The \lean{factor}s and the binomial expansion:
    \href{\leandocs\#Ring.factor}{\lean{Ring.factor}} \quad and \quad
    \href{\leandocs\#Ring.pow_add_eq_sum_factor}{\lean{Ring.pow\_add\_eq\_sum\_factor}}
  \item The vanishing of a sufficiently large power of a sum:
    \href{\leandocs\#Ring.pow_add_eq_zero_of_le}{\lean{Ring.pow\_add\_eq\_zero\_of\_le}}
  \item The first assertion of \cref{thm:fitting}:
    \href{\leandocs\#Subgroup.isNilpotent_sup_of_normal}{\lean{Subgroup.isNilpotent\_sup\_of\_normal}}
  \item The bound on the nilpotency class in \cref{thm:fitting}:
    \href{\leandocs\#Subgroup.nilpotencyClass_sup_le_of_normal}{\lean{Subgroup.nilpotencyClass\_sup\_le\_of\_normal}}
\end{itemize}

\section{Conclusion}
Viewing the collection of normal subgroups of a given group~$G$ as a semiring simplifies Fitting's theorem conceptually and formally.
Lean ensures that, even with the limited number of available axioms, the resulting semiring has a rich structure.

We found that the abstraction using a ring structure really helped both our intuition and the formalization process.

At the intuitive level, the existence of an expansion formula for the power $(H + K) ^ n$ of a sum is completely clear.
Once we identified the individual summands in the expansion, it is relatively straightforward to see that multiplication/commutator ``shrinks''.
This means that every monomial in the expansion of $(H + K) ^ n$ is a subgroup contained in both $H ^ a$ and $K ^ b$, for some $a, b$ with $a + b = n$.
From here, it is immediate to conclude that a sufficiently large power of $H + K$ is trivial, as long as sufficiently large powers of~$H$ and of~$K$ are trivial.
In fact, with a little bookkeeping, it is also possible to obtain the explicit bound mentioned in \cref{thm:fitting}.

At the formalization level, replacing the type of subgroups by a ``generic'' ring was also incredibly beneficial.
While proving the analogue of the binomial formula for non-associative rings, we could leverage \lean{Mathlib}'s extensive algebraic hierarchy to deal with the symbol manipulations and big operator formalism.
In the second part of the abstract development, we could exploit the various mixins for interactions between ring operations and their monotonicity with respect to an underlying order structure.
This again allowed us to reuse the machinery, this time the ordered hierarchy, already present in \lean{Mathlib}.

This interaction fostered what we call {\emph{type-correct creativity}}.

\bigskip
\noindent\textbf{Acknowledgements.}
I thank the \lean{Mathlib} community for infrastructure and discussions.
I tossed around some of these ideas in conversations with Thomas Browning, Inna Capdeboscq, Jireh Loreaux, Dmitriy Rumynin and Gareth Tracey.

\bigskip
\noindent\textbf{AI usage.}
Aristotle~\cite{Aristotle2025} was used to find alternative proofs.
Some of that code was reused in the final version.

Claude Opus~5~\cite{ClaudeOpus5} was used during revision, to help streamline the text and the Lean code, and to add the finishing touches.

\printbibliography

@article{Fitting1938,
  author = {Fitting, Hans},
  title = {Beitr{\"a}ge zur Theorie der Gruppen endlicher Ordnung},
  journal = {Jahresbericht der Deutschen Mathematiker-Vereinigung},
  year = {1938},
  volume = {48},
  pages  = {77--141}
}

@book{Isaacs,
  author    = {Isaacs, I. Martin},
  title     = {Finite Group Theory},
  series    = {Graduate Studies in Mathematics},
  volume    = {92},
  publisher = {American Mathematical Society},
  year      = {2008},
  address   = {Providence, RI}
}

@inproceedings{mathlib2020,
  author    = {{The mathlib Community}},
  title     = {The {L}ean {M}athematical {L}ibrary},
  booktitle = {Proceedings of the 9th {ACM} {SIGPLAN} International Conference
               on Certified Programs and Proofs},
  series    = {CPP 2020},
  publisher = {ACM},
  address   = {New Orleans, LA, USA},
  year      = {2020},
  month     = jan,
  pages     = {367--381},
  doi       = {10.1145/3372885.3373824},
  label     = {Mathlib}
}

@book{Hall1959,
  author = {Hall, Marshall},
  title = {The Theory of Groups},
  publisher = {Macmillan},
  year = {1959}
}

@book{Robinson1996,
  author = {Robinson, Derek J. S.},
  title = {A Course in the Theory of Groups},
  edition = {2},
  publisher = {Springer},
  series = {Graduate Texts in Mathematics},
  volume = {80},
  year = {1996}
}

@inproceedings{lean4,
  author    = {de Moura, Leonardo and Ullrich, Sebastian},
  title     = {The Lean 4 Theorem Prover and Programming Language (System Description)},
  booktitle = {Automated Deduction -- CADE 28},
  editor    = {Platzer, Andr\'{e} and Sutcliffe, Geoff},
  series    = {Lecture Notes in Computer Science},
  volume    = {12699},
  pages     = {625--635},
  publisher = {Springer},
  year      = {2021},
  doi       = {10.1007/978-3-030-79876-5_37},
  label     = {Lean}
}

@misc{Aristotle2025,
  author    = {Achim, Tudor and Best, Alex and Bietti, Alberto and Der, Kevin and
               F{\'e}d{\'e}rico, Math{\"i}s and Gukov, Sergei and Halpern-Leistner, Daniel and
               Henningsgard, Kirsten and Kudryashov, Yury and Meiburg, Alexander and
               Michelsen, Martin and Patterson, Riley and Rodriguez, Eric and Scharff, Laura and
               Shanker, Vikram and Sicca, Vladmir and Sowrirajan, Hari and Swope, Aidan and
               Tamas, Matyas and Tenev, Vlad and Thomm, Jonathan and Williams, Harold and
               Wu, Lawrence},
  title     = {Aristotle: {IMO}-level Automated Theorem Proving},
  year      = {2025},
  eprint    = {2510.01346},
  eprinttype = {arXiv},
  eprintclass = {cs.AI},
  label     = {Aristotle}
}

@online{ClaudeOpus5,
  author  = {{Anthropic}},
  title   = {{Claude} {Opus}~5 System Card},
  year    = {2026},
  month   = jul,
  url     = {https://www.anthropic.com/document/claude-opus-5-system-card},
  urldate = {2026-07-30},
  label   = {Claude}
}

\end{document}